\documentclass[conference, a4paper]{IEEEtran}
\usepackage{amsmath,amssymb,amsfonts}
\usepackage{algorithmic}
\usepackage{algorithm}
\usepackage{graphicx}
\usepackage{textcomp}
\usepackage{xcolor}
\usepackage{booktabs}
\usepackage{multirow}
\usepackage{url}
\usepackage{caption}
\usepackage{tikz}
\usepackage{threeparttable}
\usetikzlibrary{arrows.meta,positioning,shapes.geometric,fit}
\definecolor{roarfill}{RGB}{222,238,250}   
\definecolor{roaredge}{RGB}{ 82,142,190}   
\definecolor{roardark}{RGB}{ 30, 72,110}   
\definecolor{roarseed}{RGB}{240,247,253}   

\def\BibTeX{{\rm B\kern-.05em{\sc i\kern-.025em b}\kern-.08em
    T\kern-.1667em\lower.7ex\hbox{E}\kern-.125emX}}

\begin{document}
\bstctlcite{IEEEexample:BSTcontrol}
\title{\textbf{ROAR}: \textbf{R}estricted 2-\textbf{O}pt \textbf{A}nd \textbf{R}uin-and-recreate Heuristic for Large and Extreme Scale Traveling Repairman Problem}

\author{
\IEEEauthorblockN{
 Khoa Hoang Dang Nguyen\textsuperscript{1},\quad
 Nam Khoa Pham\textsuperscript{1},\quad
 Khanh Dinh Nam Nguyen\textsuperscript{1},\quad
 Bac Ba Thanh Nguyen\textsuperscript{1,2}}
\IEEEauthorblockA{
\textsuperscript{1}\textit{College of Engineering and Computer Science}, \textsuperscript{2}\textit{Center for AI Research}, \textit{VinUniversity}, Hanoi, Vietnam\\
\{25khoa.nhd, 25nam.pk, 25khanh.ndn, 25bac.nbt\}@vinuni.edu.vn}
}

\maketitle

\begin{abstract}
The Traveling Repairman Problem minimizes cumulative arrival times rather than the total route length, as in the classical Traveling Salesman Problem. The literature contains several legacy metaheuristics for the TRP; however, most are evaluated on small- and medium-scale instances and show limited scalability for large and extreme datasets. GILS-RVND, the current gold standard in metaheuristic algorithms, performs well on small and medium instances but struggles to scale due to the exhaustive neighborhood search in the RVND phase. To address this limitation, we introduce \textbf{R}estricted 2-\textbf{O}pt \textbf{A}nd \textbf{R}uin-and-recreate \textbf{(ROAR)}, a Large Neighborhood Search–based method that integrates ruin-and-recreate with greedy insertion and a restricted, insertion-anchored 2-opt operator. Computational experiments demonstrate that ROAR produces high-quality solutions for large-scale instances and outperforms existing methods under strict time limits. We further introduce an $\mathcal{O}(1)$-time, $\mathcal{O}(1)$-space algorithm for evaluating 2-opt moves during a contiguous systematic sweep, which, to the best of our knowledge, is unprecedented. Finally, we evaluate ROAR on instances with up to $200{,}000$ nodes, substantially exceeding the largest instances previously reported in the literature ($11{,}849$ nodes).
\end{abstract}

\section{Introduction}
\label{sec:intro}
The Traveling Repairman Problem (TRP), equivalently the Minimum Latency Problem (MLP), aims to find a route visiting every vertex that minimizes the sum of arrival times rather than, as in the Traveling Salesman Problem (TSP), the Hamiltonian cycle length. The objective is cumulative: any local reordering shifts every subsequent arrival, so cost responds non-locally. It models problems that prioritize waiting time over tour length. Some typical situations are post-disaster rescues or customer-centric delivery systems.

The problem is NP-hard; exact solvers are generally applicable up to about $200$ nodes \cite{bulhoes2018}; larger graphs demand heuristics. Existing metaheuristics are assessed almost exclusively on small and medium instances, the largest reaching $11{,}849$ nodes \cite{yelmewad2021}. GILS-RVND \cite{silva2012}, our principal benchmark, is constrained both by the exhaustive RVND phase and by $\mathcal{O}(N^2)$ auxiliary data structures that become prohibitive in terms of memory consumption at larger problem scales (Section~\ref{sec:related}). PDLSH \cite{yelmewad2021} attains substantially greater problem sizes via GPU acceleration, albeit at the cost of reduced solution quality and the requirement of a dedicated hardware accelerator.

In this work, we propose a lightweight, CPU-native heuristic returning strong solutions within a fixed, predefined time budget at extreme scale. It revisits \textbf{Incremental Local Optimization} (ILO) \cite{caseau1999}: interleaving optimization with insertion rather than deferring it to a post-optimization pass. In short, we propose:
\begin{enumerate}
    \item \textbf{R}estricted 2-\textbf{O}pt \textbf{A}nd \textbf{R}uin-and-recreate (\textbf{ROAR}), embedding a \emph{restricted 2-opt} inside the recreation phase of an LNS scheme.
    \item A 2-opt evaluation algorithm attaining $\mathcal{O}(1)$ time-complexity per move during a contiguous systematic sweep with only $\mathcal{O}(1)$ auxiliary space;
    \item The first metaheuristic baseline for extreme-scale TRP, up to $200{,}000$ nodes; To the best of our knowledge, the largest instance experimented to date comprises \(11{,}849\) nodes \cite{yelmewad2021}.
\end{enumerate}

\section{Related Work}
\label{sec:related}

GILS-RVND \cite{silva2012}, which couples a Greedy Randomized Adaptive Search Procedure with Iterated Local Search and Randomized Variable Neighborhood Descent (RVND) and iteratively searches each neighborhood until no improvement is found. Despite being unmatched on small and medium instances, it fails on two fronts: exhaustive descent grows prohibitive as $N$ rises, and a memory wall, since a TRP 2-opt move costs $\mathcal{O}(N)$ naively and its $\mathcal{O}(1)$ 2-opt evaluator relies on precomputed $\mathcal{O}(N^2)$ matrices. Subsequent methods have sought to address these limitations. DBMEA \cite{tuuszabo2020} pairs a discrete bacterial memetic evolutionary scheme with 2-opt and 3-opt, reports new best-known values at $500$ and $1{,}000$ nodes, and cuts auxiliary space for the constant-time 2-opt evaluator to $\mathcal{O}(N)$ via cumulative sums. MDM-GILS-RVND \cite{santana2022} adds data mining, matching the baseline on small instances and improving marginally to $1{,}379$ nodes. PDLSH \cite{yelmewad2021} parallelizes move evaluation across GPU threads, scaling beyond $10{,}000$ nodes, but sacrifices solution quality.

More recently, NeuFACO \cite{11365075} introduces a neural-guided Ant Colony Optimization (ACO) framework, combining learned heuristic guidance with ACO and local search to improve solution quality and scalability for TSP problem. DyNACO \cite{tran2026beyond} introduces dynamic neural guidance that adapts to the evolving pheromone distribution and incumbent solution during ACO search, enabling optimization of large-scale problems.

Large Neighborhood Search (LNS) \cite{shaw1998}, also known as ruin-and-recreate \cite{schrimpf2000}, and its adaptive variant (ALNS) \cite{ropke2006} are fundamentals in routing problems: each iteration \emph{ruins} part of an incumbent and \emph{recreates} it by a repair operator, escaping optima beyond move-based search. Within this paradigm, 2-opt local search \cite{croes1958} is commonly applied as a post-optimization step after repair to further improve solution quality \cite{voigt2025}.

The interleaving of local search procedures with insertion operations is less common. GENI \cite{gendreau1992} restructures the tour as each vertex enters; its US post-optimization applies restricted local search after one vertex is removed and re-inserted; Caseau and Laburthe \cite{caseau1999} formalized this as ILO, a bounded local search run immediately after each insertion on the partial solution. Adoption remains mostly confined to VRP and TSP, and to construction rather than the recreate step of an LNS loop; only \cite{caseau2001} couples the two for the VRP. To our knowledge, ILO has not been applied to the TRP.

\section{Problem Formalization}
TRP is formally defined on a symmetric, complete, weighted graph $G = (V, E)$, where $V$ is the set of $N$ vertices and $E$ is the set of connecting edges. Vertex $0$ designates the starting depot. Let $c(u, v)$ represent the cost to travel between vertex $u$ and vertex $v$. A valid routing sequence is a permutation of the vertices $\pi = (\pi_0, \pi_1, \dots, \pi_{N-1})$, constrained such that the sequence originates at the depot ($\pi_0 = 0$). The arrival time, or latency, at the $k$-th node ($k$ starts from 0) in the route is the cumulative travel time from the depot to that node, defined as:
\begin{equation}
    L_k = \sum_{i=0}^{k-1} c(\pi_i, \pi_{i+1}).
\end{equation}

In the literature, the TRP objective function has two variants: \textbf{Depot-Inclusive}, which includes depot latency in the final return, and \textbf{Depot-Exclusive}, which does not. We adopt the \textbf{Depot-Inclusive} formulation, as in \cite{silva2012} for some datasets, our primary benchmark. Here, $Z(\pi)$ denotes the tour cost: 
\begin{equation}
    Z(\pi) = \sum_{k=1}^{N-1} L_k + L_{N-1} + c(\pi_{N-1}, \pi_{0})
\end{equation}

\begin{table}[t]
\centering
\begin{threeparttable}
\caption{TRP metaheuristics comparison}
\label{tab:landscape}
\setlength{\tabcolsep}{5pt}
\renewcommand{\arraystretch}{1.15}
\footnotesize
\begin{tabular}{@{}l r c c@{}}
\toprule
Method & Max $N$ & $\mathcal{O}(1)$ 2-opt evaluation & Aux.\ space \\
\midrule
GILS-RVND \cite{silva2012}      &   1{,}000 & \checkmark & $\mathcal{O}(N^2)$ \\
DBMEA \cite{tuuszabo2020}       &   1{,}000 & \checkmark & $\mathcal{O}(N)$   \\
MDM-GILS-RVND \cite{santana2022}&   1{,}379 & \checkmark & $\mathcal{O}(N^2)$ \\
PDLSH \cite{yelmewad2021}       &  11{,}849 & $\times$   & -- \\
\addlinespace[1pt]
\textbf{ROAR} (ours)            & 200{,}000 & \checkmark & $\mathcal{O}(1)$   \\
\bottomrule
\end{tabular}
\begin{tablenotes}[flushleft]\scriptsize
\item \emph{Max $N$}: largest experimented instance.
\emph{Aux.\ space}: extra memory for this operation.

\end{tablenotes}
\end{threeparttable}
\end{table}
\section{Methodology}

ROAR follows the ruin-and-recreate paradigm: from a seed tour it repeatedly removes a stochastic subset of nodes and reinserts them greedily, interleaving a restricted 2-opt within the repair phase to remove the crossings each insertion introduces---the ILO scheme. 
\begin{figure}[t]
\centering
\begin{tikzpicture}[
  node distance=5.8mm,
  box/.style  ={draw=roaredge, line width=0.45pt, fill=roarfill,
                rounded corners=1.6pt, font=\footnotesize, align=center,
                text=roardark, inner xsep=3.5pt, inner ysep=2.4pt},
  seed/.style ={box, fill=roarseed},
  flow/.style ={-{Latex[length=1.7mm,width=1.4mm]}, semithick, draw=roardark},
  lbl/.style  ={font=\scriptsize, text=roardark, inner sep=1.5pt}
]
\node[seed]             (s) {Nearest Neighbour\\[-1pt]{\scriptsize\color{roaredge}seed}};
\node[box, right=of s]  (r) {Ruin\\[-1pt]{\scriptsize\color{roaredge}random}};
\node[box, right=of r]  (i) {Insert\\[-1pt]{\scriptsize\color{roaredge}greedy}};
\node[box, right=of i]  (t) {2-opt\\[-1pt]{\scriptsize\color{roaredge}restricted}};
\node[box, right=of t]  (a) {Accept\\[-1pt]{\scriptsize\color{roaredge}greedy}};

\draw[flow] (s) -- (r);  \draw[flow] (r) -- (i);
\draw[flow] (i) -- (t);  \draw[flow] (t) -- (a);
\draw[flow] (a) -- ++(7.5mm,0) coordinate (e);
\node[lbl, right=0.5pt of e] {$\pi^\star$};

\draw[flow] (t.north) to[out=118,in=62]
      node[lbl, midway, above=-1pt] (lp) {$v\in R$} (i.north);
\node[draw=roaredge, dashed, dash pattern=on 1.2pt off 1.2pt,
      rounded corners=2pt, inner sep=2pt, fit=(i)(t)(lp),
      label={[lbl, roaredge, fill=white, inner sep=1.5pt]above:ILO}] {};

\path (r.south) -- (a.south) coordinate[midway] (m);
\draw[flow] (a.south) -- ++(0,-4.4mm) -| (r.south);
\node[lbl, below=4.2mm of m, inner sep=1.5pt] {repeat};
\end{tikzpicture}
\caption{The ROAR pipeline}
\label{fig:roar}
\end{figure}
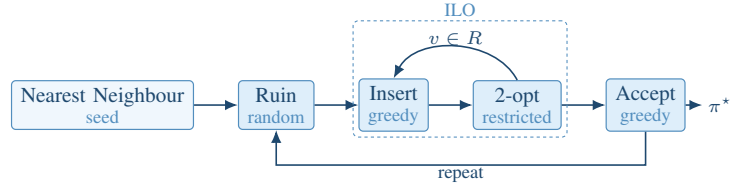
\subsection{Initialization}
Any greedy construction may serve; we adopt the classical nearest-neighbor heuristic \cite{rosenkrantz1977}, appending from the depot the nearest unvisited node until the initial permutation $\pi$ is complete.

\subsection{Search Loop and Termination}
The loop runs until either of two conditions is met. The first is convergence: the best-found cost is tracked over epochs of length $E(N)$, and the search halts after $k$ consecutive non-improving epochs. The second is a wall-clock budget $T(N)$, which encodes the operational premise of this work---a usable solution within a predictable timeframe, irrespective of convergence. Both scale with $N$ (Eqs.~\eqref{eq:epoch},~\eqref{eq:budget}); the incumbent is returned as soon as either fires.

\subsection{The ROAR Cycle}
\label{sec:cycle}
Algorithm~\ref{alg:roar} summarizes a ROAR cycle. Each cycle draws a removal count $N_{\text{remove}} \sim \mathcal{U}\{1,\dots,N_{\text{remove}}^{\max}\}$, where
\begin{equation}
    N_{\text{remove}}^{\max} = \min\big(\lfloor \alpha_{\text{destroy}} \cdot (N_{\text{total}} - 1) \rfloor,\; N_{\text{total}} - 2,\; C_{\text{limit}} \big).
\end{equation}
Here $\alpha_{\text{destroy}} \in (0,1)$ controls the aggressiveness of the ruin and $C_{\text{limit}}$ caps the per-cycle runtime at extreme scale; the term $N_{\text{total}} - 2$ retains at least the depot and one anchor node.

\subsubsection{Ruin (Destroy)} The $N_{\text{remove}}$ nodes are drawn uniformly from the non-depot positions and excised, leaving a \emph{partial route} $\pi'$ and a pool $R$ of removed nodes.

\subsubsection{Recreate (Repair)} The nodes of $R$ are reinserted one at a time in random order (arbitrary insertion \cite{rosenkrantz1977}). To insert $v$ between the consecutive nodes $\pi_{i-1}$ and $\pi_i$, let $T_{i-1}$ be the cumulative travel time from the depot to $\pi_{i-1}$ and let
\begin{equation}
    \delta_i = c(\pi_{i-1}, v) + c(v, \pi_i) - c(\pi_{i-1}, \pi_i)
\end{equation}
be the induced detour. The objective then increases by
\begin{equation}
\label{eq:insertion}
    \Delta L(v, i) = \underbrace{\big[\, T_{i-1} + c(\pi_{i-1}, v) \,\big]}_{\text{arrival time of } v} \;+\; \underbrace{m_i \cdot \delta_i}_{\text{cascaded delay}},
\end{equation}
where $m_i$ counts the nodes after the insertion slot, including the depot return under the depot-inclusive objective. The node is placed at $i^\star = \arg\min_i \Delta L(v, i)$. As $T_{i-1}$ is carried incrementally, each candidate costs $\mathcal{O}(1)$ and a full scan $\mathcal{O}(|\pi'|)$. The detailed algorithm for this is shown in Algorithm ~\ref{alg:insert}.
\begin{algorithm}[t]
\caption{Insertion of node $v$ into $\pi'$}
\label{alg:insert}
\begin{algorithmic}[1]
\STATE \textbf{Input:} partial route $\pi'=(\pi'_0,\dots,\pi'_{m-1})$ with depot $\pi'_0$, node $v$
\STATE $T \leftarrow 0$,\ \ $\Delta^\star \leftarrow \infty$,\ \ $i^\star \leftarrow 1$ \quad // $T$: prefix time to slot
\FOR{$i = 1$ to $m$}
    \STATE $a \leftarrow \pi'_{i-1}$,\ \ $b \leftarrow \pi'_{\,i \bmod m}$ \quad // $b=$ depot when $i=m$
    \STATE $\delta \leftarrow c(a,v)+c(v,b)-c(a,b)$
    \STATE $\Delta \leftarrow T + c(a,v) + (m-i+1)\,\delta$ \quad // Eq.~\eqref{eq:insertion}
    \IF{$\Delta < \Delta^\star$}
        \STATE $\Delta^\star \leftarrow \Delta$,\ \ $i^\star \leftarrow i$
    \ENDIF
    \STATE $T \leftarrow T + c(a,b)$ \quad // advance prefix time
\ENDFOR
\STATE insert $v$ into $\pi'$ at position $i^\star$
\end{algorithmic}
\end{algorithm}
\subsubsection{Restricted 2-opt} 
\begin{figure}[t]
\centering
\begin{tikzpicture}[x=6.2mm, y=1mm,
  dot/.style  ={circle, draw=roaredge, fill=roarfill, line width=0.35pt,
                inner sep=0pt, minimum size=1.7mm},
  vdot/.style ={circle, draw=roardark, fill=roardark, inner sep=0pt, minimum size=2.3mm},
  dep/.style  ={rectangle, draw=roardark, fill=white, line width=0.5pt,
                inner sep=0pt, minimum size=1.9mm},
  tour/.style ={draw=roaredge!55, line width=0.5pt},
  anch/.style ={draw=roardark, line width=1.2pt},
  sweep/.style={-{Latex[length=1.5mm,width=1.2mm]}, draw=roaredge, line width=0.6pt},
  tick/.style ={draw=roaredge, line width=0.4pt},
  drop/.style ={draw=roardark, line width=0.4pt},
  lbl/.style  ={font=\tiny, text=roardark, inner sep=1pt}
]
\foreach \k in {1,2,3,4,5,7,8,9,10,11}{\node[dot] (n\k) at (\k,0) {};}
\node[dep] (n0) at (0,0) {};  \node[dep] (n12) at (12,0) {};
\node[vdot] (n6) at (6,0) {};
\foreach \k [evaluate=\k as \j using int(\k+1)] in {0,...,11}{\draw[tour] (n\k) -- (n\j);}
\draw[anch] (n5) -- (n6);  \draw[anch] (n6) -- (n7);
\node[lbl] at (6.42,-2.0) {$v$};
\node[lbl] at (0,-2.2) {$0$};  \node[lbl] at (12,-2.2) {$0$};
\foreach \k in {1,...,4} {\draw[tick] (\k,3.2) -- (\k,4.8);}
\foreach \k in {7,...,11}{\draw[tick] (\k,3.2) -- (\k,4.8);}
\draw[sweep] (5,4) -- (0.15,4);   \node[lbl] at (2.7,6.0) {$F_1$};
\draw[sweep] (6,4) -- (11.85,4);  \node[lbl] at (9.2,6.0) {$F_2$};
\draw[drop] (5,1.0) -- (5,4);  \draw[drop] (6,1.2) -- (6,4);
\node[lbl] at (5.5,6.0) {$e_{in}$};
\foreach \k in {1,...,5} {\draw[tick] (\k,-5.2) -- (\k,-6.8);}
\foreach \k in {8,...,11}{\draw[tick] (\k,-5.2) -- (\k,-6.8);}
\draw[sweep] (6,-6) -- (0.15,-6);   \node[lbl] at (3.0,-8.2) {$F_3$};
\draw[sweep] (7,-6) -- (11.85,-6);  \node[lbl] at (9.6,-8.2) {$F_4$};
\draw[drop] (6,-1.2) -- (6,-6);  \draw[drop] (7,-1.0) -- (7,-6);
\node[lbl] at (6.5,-8.2) {$e_{out}$};
\end{tikzpicture}
\caption{Restricted 2-opt at an insertion: one cut edge is fixed at the new node $v$—its entry edge $e_{in}=(\pi_{i-1},v)$ or exit edge $e_{out}=(v,\pi_i)$—while the other moves outward one node at a time, forming the four contiguous families $F_1$–$F_4$.}
\label{fig:r2opt}
\end{figure}
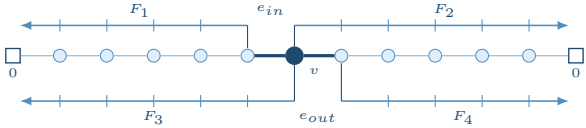

After that, ROAR runs a 2-opt search restricted to reversals in which one broken edge is incident to the inserted node, rather than the full $\mathcal{O}(N^2)$ neighborhood. Anchoring on the entry and exit edges of the insertion point yields four contiguous families of reversals; within each, the second cut edge advances monotonically along the remainder of the route, so the segment latency terms follow a recurrence and every candidate is assessed in $\mathcal{O}(1)$ (Section~\ref{sec:o1}). The best improving reversal across the four families, if any, is applied (Fig.~\ref{fig:r2opt}).
\subsubsection{Acceptance policy} The new tour is accepted only if it yields improvement over the pre-ruin tour.

\begin{algorithm}[t]
\caption{ROAR: Ruin-and-Recreate with Restricted 2-opt}
\label{alg:roar}
\begin{algorithmic}[1]
\STATE \textbf{Input:} seed tour $\pi$, destruction level $\alpha_{\text{destroy}}$, removal cap $C_{\text{limit}}$
\STATE $Z \leftarrow Z(\pi)$ \quad // track the current cost
\STATE $N_{\text{remove}}^{\max} \leftarrow \min\!\big(\lfloor \alpha_{\text{destroy}}(N_{\text{total}}{-}1)\rfloor,\; N_{\text{total}}{-}2,\; C_{\text{limit}}\big)$
\WHILE{stopping criteria not met}
    \STATE draw $N_{\text{remove}} \sim \mathcal{U}\{1,\dots,N_{\text{remove}}^{\max}\}$
    \STATE $(\pi', R) \leftarrow \textsc{Ruin}(\pi, N_{\text{remove}})$ \quad // random removals
        
    \FOR{each $v \in R$ in random order}     
        \STATE insert $v$ at $\arg\min_i \Delta L(v,i)$ \quad // Eq.~\eqref{eq:insertion}
        \STATE apply best anchored restricted 2-opt at $v$ \quad // Eq.~\eqref{eq:deltaZ}
    \ENDFOR
    \IF{$Z(\pi') < Z$}
        \STATE $\pi \leftarrow \pi'$,\ \ $Z \leftarrow Z(\pi')$ \quad // greedy accept
    \ENDIF
\ENDWHILE
\STATE \textbf{return} $\pi$
\end{algorithmic}
\end{algorithm}
\subsection{Constant-time 2-opt Evaluation}
\label{sec:o1}
A naive 2-opt evaluation on the TRP costs $\mathcal{O}(N)$. Existing works trade space for time.
GILS-RVND precomputes $\mathcal{O}(N^{2})$ matrices for constant-time evaluation; DBMEA attains the same
cost in $\mathcal{O}(N)$ space via cumulative sums; the concatenation framework of
\cite{vidal2014} composes duration, cumulative cost, and delay weight inductively across subsequences.
We further optimize the memory consumption in this phase. The restricted 2-opt never evaluates an isolated move: it anchors one
fracture edge and advances the other monotonically. With this setting, we can push the auxiliary space complexity to constant.

This algorithm requires the second edge to advance monotonically. Iterating the anchor over all $N$ edges thus covers the full 2-opt neighborhood in $\mathcal{O}(N^{2})$ time and $\mathcal{O}(1)$ space. Thus, this algorithm can be integrated into methods that use a full 2-opt neighborhood scan. This is applied in the LNS+2-opt algorithm in the experiments.

\subsubsection{Mathematical Decomposition of Latency Change}
Let $\pi = (\pi_0, \pi_1, \dots, \pi_{N-1})$ be a valid solution. Given that $i < j < N$, a standard 2-opt move removes two edges, $e_1 = (\pi_i, \pi_{i+1})$ and $e_2 = (\pi_j, \pi_{j+1})$, and reconnects them with $e'_1 = (\pi_i, \pi_j)$ and $e'_2 = (\pi_{i+1}, \pi_{j+1})$, thereby reversing the middle segment $P_{mid} = (\pi_{i+1}, \dots, \pi_j)$ into $P'_{mid} = (\pi_j, \dots, \pi_{i+1})$. The total change in latency, $\Delta Z = Z(\pi') - Z(\pi)$, partitions into three independent components:
\begin{enumerate}
    \item \textbf{Prefix stability:} The arrival times of all nodes in the prefix $(\pi_0, \dots, \pi_i)$ remain entirely unchanged.
    \item \textbf{Suffix shift:} At the entry of the reversed segment, the edge length changes by $\delta_{entry} = c(\pi_i, \pi_j) - c(\pi_i, \pi_{i+1})$. Because latency is cumulative, this shift cascades linearly across all $K_b$ nodes evaluated after index $i$. The new exit edge introduces a secondary shift $\delta_{exit} = c(\pi_{i+1}, \pi_{j+1}) - c(\pi_j, \pi_{j+1})$, which cascades across all $K_c$ nodes appearing after index $j$. 
    \item \textbf{Internal reversal:} The latency accumulated within the middle segment is restructured. Let $TRP_{first}$ denote the relative latency of $P_{mid}$ in its original orientation, and $TRP_{rev}$ the relative latency of $P'_{mid}$ when reversed.
\end{enumerate}
Synthesizing these components yields the exact cost delta:
\begin{equation}
\label{eq:deltaZ}
\Delta Z = \delta_{entry}K_b + \delta_{exit}K_c + TRP_{rev} - TRP_{first}.
\end{equation}
Here, $K_b$ and $K_c$ are inferred: $K_b = N - i$ and $K_c = N - j$. Note that this formula incorporates the return-to-depot latency.

\subsubsection{Continuous Sweep Algorithm}
Evaluating Eq.~\eqref{eq:deltaZ} for an isolated pair $(i,j)$ still requires an $\mathcal{O}(N)$ traversal to compute $TRP_{rev}$ and $TRP_{first}$. We achieve $\mathcal{O}(1)$ time by structuring the search as a systematic contiguous sweep: anchoring the first fracture edge at $i$ and advancing the second fracture edge $j \to j+1$ step by step, the segment terms are aggregated incrementally. Letting $L_{tsp}$ track the length of the reversed segment, from node $i+1$ to $j$ (endpoint inclusive). As $j-1$ increments to $j$, we update these terms as follows, \textbf{with the variables update order as shown}:
\begin{align}
    L_{tsp} &\leftarrow L_{tsp} + c(\pi_{j-1}, \pi_j), \\
    TRP_{rev} &\leftarrow TRP_{rev} + c(\pi_{j-1}, \pi_j) \times (j - i - 1), \\
    TRP_{first} &\leftarrow TRP_{first} + L_{tsp}.
\end{align}
In the trivial case $j = i+1$, all terms are 0, so computing $\Delta Z$ costs $\mathcal{O}(1)$ per movement, assuming we evaluate a full chain of movements where $i$ is fixed and $j$ increments by 1 in each step. The algorithm where $j$ is fixed while $i$ decrements can be derived from this base algorithm with minor modifications.
\subsubsection{Limitations}
This algorithm holds only under a contiguous sweep: one fracture edge fixed, the
other advancing monotonically. Arbitrary 2-opt moves cannot be evaluated in
constant time by this method. ROAR's restricted 2-opt satisfies this access pattern by construction.
\section{Computational Experiments}

\subsection{Environment and Baselines}
ROAR is benchmarked against \textbf{GILS-RVND}, the standard CPU
metaheuristic for the TRP, and \textbf{PDLSH}, a GPU-parallel local
search; in addition, we also benchmark against a standard combination of LNS and 2-opt, using 2-opt as a post-repair operator.

All runs use a single workstation with an Intel Core i7-13620H CPU, $32$\,GB DDR5,
and an NVIDIA RTX~4060 GPU ($8$\,GB VRAM) employed by PDLSH alone; every reported
cost is a mean over $10$ independent runs.

Our GILS-RVND implementation largely adapts the
codebase\footnote{\url{https://github.com/franciscunha/MLP}} with some modifications. We replace the GRASP-based initialization of GILS-RVND with a nearest-neighbor initialization because, for large-scale instances, GRASP adds substantial computational overhead and unfairly degrades performance under limited runtime. This does not degrade the baseline: GILS-RVND still
reproduces the published optimum on every instance up to $150$ nodes, the sole
exception being \texttt{pr107}, within $0.07\%$ on average. PDLSH is the authors'
official implementation, aligned as described in Section~\ref{sec:pdlsh-fix}.
Latencies are rounded to five significant digits; full values, re-implementations, ROAR cycle counts (for probabilistic reproducibility)
and our codebase are available
online\footnote{\url{https://github.com/khoa-bot/ROARimplementation}}. 
\subsection{Instances and Reference Values}
Instances are drawn from TSPLIB~\cite{reinelt1991,reinelt1995}, the VLSI
collection~\cite{rohe_vlsi}, and the TSP-art set~\cite{kaplan2005,bosch_art}, and are
grouped into a \emph{regular} set of $52$--$11{,}849$ nodes, a \emph{large} set of
$13{,}509$--$33{,}810$ nodes, and an \emph{extreme} set of $32{,}892$--$200{,}000$
nodes. Reference values, either best-known or proven optimal, are available only for the regular set, as this is the first time instances from large and extreme sets are being experimented with. Our metric assigns edge weights according to each instance’s \textbf{EDGE\_WEIGHT\_TYPE} specification.

\subsection{Objective Alignment}
\label{sec:pdlsh-fix}
The official PDLSH implementation departs from our benchmark metric—adopted by \cite{silva2012} (for the Hamiltonian Circuit Variant) and \cite{abeledo2013}—in two key respects. First, its Euclidean routine assigns the result of \texttt{sqrtf} to a \texttt{long} integer, truncating rather than rounding fractional edge lengths for \texttt{EUC\_2D} instances. Second, its objective calculation accumulates arrival times only at customer nodes, omitting the final return to the depot. A mistake about mixing these objectives is apparent in \cite{yelmewad2021}: because both deviations systematically underestimate total latency, comparing the results of \cite{yelmewad2021} directly against \cite{abeledo2013} and \cite{silva2012} leads to invalid claims of ``new best-known'' solutions on instances proven optimal, such as \texttt{berlin52}, in \cite{yelmewad2021}. For proper comparison, we modify the edge-weight computation and the latency metric to align with our benchmark.

However, because this is the sole work that reports experiments on instances comprising up to more than $11{,}000$ nodes, we keep the costs for the corresponding datasets reported in \cite{yelmewad2021} as an optimistic reference baseline in Table~\ref{tab:main}.

\subsection{Hyperparameters and Time Budget}
ROAR and LNS+2-opt use one fixed configuration throughout: destruction level
$\alpha_{\text{destroy}} = 0.2$, removal cap $C_{\text{limit}} = 300$, and a patience
of $10$ consecutive non-improving epochs, an \emph{epoch} being a wall-clock interval
of $E(N)$ seconds. The epoch length scales with instance size,
\begin{equation}
\label{eq:epoch}
E(N)=
\begin{cases}
1 + 29\cdot\dfrac{N-1000}{99{,}000}, &  N \le 100{,}000,\\[6pt]
30, & N > 100{,}000,
\end{cases}
\end{equation}
and a wall-clock budget $T(N)$ bounds every method:
\begin{equation}
\label{eq:budget}
T(N)=
\begin{cases}
100, & N < 1000,\\[3pt]
100 + \dfrac{N-1000}{110}, & 1000 \le N \le 100{,}000,\\[6pt]
1000, & N > 100{,}000.
\end{cases}
\end{equation}

\subsection{Regular Set}
\label{sec:res-regular}

\begin{table*}[t]
\centering
\begin{threeparttable}
\caption{Regular set ($52$--$11{,}849$ nodes) and large set ($13{,}509$--$33{,}810$
nodes). Costs are divided by the \emph{Scale} factor; $t$: mean wall-clock time (s)
over $10$ runs. Lowest cost per instance in \textbf{bold}; \underline{underlined}
entries fall below the reference value. LNS+2-opt applies the 2-opt scan after each
repair rather than at each insertion. The last three columns give
$100\,(Z_{\text{ROAR}}-Z_{\text{base}})/Z_{\text{base}}$, so a negative entry
favours ROAR; in the \emph{vs.\ BKS} column a negative (\textbf{bold}) entry marks a
cost lower than reference value. Reference values exist only for the regular set.
\textsc{oom}: out of memory.}
\label{tab:main}
\renewcommand{\arraystretch}{1.15}
\setlength{\tabcolsep}{3pt}
\footnotesize
\begin{tabular}{@{}l r c r r r r r r r r r r r r@{}}
\toprule
 & & & & \multicolumn{2}{c}{GILS-RVND \cite{silva2012}} & \multicolumn{2}{c}{PDLSH \cite{yelmewad2021}} & \multicolumn{2}{c}{LNS+2-opt} & \multicolumn{2}{c}{ROAR (ours)} & \multicolumn{3}{c}{ROAR gap (\%)} \\
\cmidrule(lr){5-6}\cmidrule(lr){7-8}\cmidrule(lr){9-10}\cmidrule(lr){11-12}\cmidrule(lr){13-15}
Instance & $N$ & Scale & Best-known & Cost & $t$(s) & Cost & $t$(s) & Cost & $t$(s) & Cost & $t$(s) & \shortstack[r]{vs.\\BKS} & \shortstack[r]{vs.\\GILS} & \shortstack[r]{vs.\\LNS+2-opt} \\
\midrule
\multicolumn{15}{@{}l}{\emph{Regular set}}\\
berlin52 & 52 & $10^{5}$ & 1.4372\textsuperscript{\textdagger} & \textbf{1.4372} & 0.1 & 1.4833 & 0.0 & 1.4641 & 11.5 & 1.4565 & 8.0 & 1.34 & 1.34 & $-0.52$ \\
st70 & 70 & $10^{4}$ & 2.0557\textsuperscript{\textdagger} & \textbf{2.0557} & 0.3 & 2.2083 & 0.0 & 2.1655 & 8.4 & 2.0698 & 8.7 & 0.69 & 0.69 & $-4.42$ \\
kroA100 & 100 & $10^{6}$ & 0.98313\textsuperscript{\textdagger} & \textbf{0.98313} & 1.3 & 1.0256 & 0.0 & 1.0196 & 10.4 & 1.0174 & 10.8 & 3.49 & 3.49 & $-0.22$ \\
pr107 & 107 & $10^{6}$ & 2.0266\textsuperscript{\textdagger} & \textbf{2.0280} & 1.2 & 2.0302 & 0.0 & 2.0298 & 8.3 & 2.0299 & 8.9 & 0.16 & 0.09 & 0.00 \\
ch150 & 150 & $10^{5}$ & 4.4442\textsuperscript{\S} & \textbf{4.4442} & 4.3 & 4.8709 & 0.0 & 4.6896 & 14.1 & 4.5795 & 12.8 & 3.04 & 3.04 & $-2.35$ \\
d198 & 198 & $10^{6}$ & 1.1860\textsuperscript{\P} & \textbf{1.1861} & 12.0 & 1.2342 & 0.0 & 1.2019 & 10.9 & 1.1998 & 11.7 & 1.16 & 1.16 & $-0.17$ \\
ts225 & 225 & $10^{7}$ & 1.3240\textsuperscript{\P} & \textbf{1.3242} & 9.3 & 1.3793 & 0.0 & 1.3499 & 8.5 & 1.3297 & 9.4 & 0.43 & 0.42 & $-1.50$ \\
a280 & 280 & $10^{5}$ & 3.4699\textsuperscript{\P} & \textbf{3.4715} & 43.5 & 3.6818 & 0.1 & 3.6363 & 18.8 & 3.5410 & 22.9 & 2.05 & 2.00 & $-2.62$ \\
pcb442 & 442 & $10^{7}$ & 1.0302\textsuperscript{\P} & \textbf{1.0332} & 100.0 & 1.0728 & 0.1 & 1.0523 & 21.0 & 1.0450 & 22.3 & 1.44 & 1.14 & $-0.69$ \\
att532 & 532 & $10^{6}$ & 5.6130\textsuperscript{\P} & \textbf{5.6439} & 100.0 & 6.2646 & 0.9 & 5.9249 & 31.1 & 5.8682 & 31.6 & 4.55 & 3.97 & $-0.96$ \\
rat783 & 783 & $10^{6}$ & 3.2759\textsuperscript{\P} & \textbf{3.3413} & 100.0 & 3.4693 & 1.9 & 3.4060 & 38.4 & 3.3945 & 46.0 & 3.62 & 1.59 & $-0.34$ \\
dsj1000ceil & 1{,}000 & $10^{9}$ & 7.6464\textsuperscript{\P} & \textbf{7.9527} & 100.0 & 8.6710 & 5.1 & 8.1666 & 62.9 & 8.0882 & 65.2 & 5.78 & 1.70 & $-0.96$ \\
u1432 & 1{,}432 & $10^{8}$ & 1.1187\textsuperscript{\textdaggerdbl} & \textbf{\underline{1.0823}} & 103.9 & 1.1207 & 19.8 & \underline{1.1003} & 96.9 & \underline{1.0966} & 92.6 & $\mathbf{-1.98}$ & 1.32 & $-0.34$ \\
d1655 & 1{,}655 & $10^{7}$ & 4.9465\textsuperscript{\textdaggerdbl} & \textbf{\underline{4.8331}} & 106.0 & 4.9715 & 20.8 & \underline{4.8912} & 67.0 & \underline{4.8504} & 92.6 & $\mathbf{-1.94}$ & 0.36 & $-0.83$ \\
pr2392 & 2{,}392 & $10^{8}$ & 4.6486\textsuperscript{\textdaggerdbl} & \underline{4.4409} & 112.7 & 4.6584 & 113.0 & \underline{4.4618} & 109.2 & \textbf{\underline{4.4394}} & 112.7 & $\mathbf{-4.50}$ & $-0.03$ & $-0.50$ \\
pcb3038 & 3{,}038 & $10^{8}$ & 2.0952\textsuperscript{\textdaggerdbl} & \textbf{\underline{2.0554}} & 118.5 & 2.1237 & 120.1 & \underline{2.0604} & 118.5 & \underline{2.0570} & 118.5 & $\mathbf{-1.82}$ & 0.08 & $-0.17$ \\
fnl4461 & 4{,}461 & $10^{8}$ & 4.0618\textsuperscript{\textdaggerdbl} & 4.1252 & 131.5 & 4.2290 & 131.8 & \underline{4.0431} & 131.5 & \textbf{\underline{4.0283}} & 131.5 & $\mathbf{-0.82}$ & $-2.35$ & $-0.37$ \\
rl5934 & 5{,}934 & $10^{9}$ & 1.6335\textsuperscript{\textdaggerdbl} & \underline{1.6306} & 144.9 & 1.6778 & 155.0 & \underline{1.5995} & 144.9 & \textbf{\underline{1.5886}} & 144.9 & $\mathbf{-2.75}$ & $-2.58$ & $-0.68$ \\
pla7397 & 7{,}397 & $10^{10}$ & 7.0429\textsuperscript{\textdaggerdbl} & 7.4927 & 158.2 & 7.9659 & 178.8 & 7.1918 & 158.2 & \textbf{\underline{6.9388}} & 158.2 & $\mathbf{-1.48}$ & $-7.39$ & $-3.52$ \\
rl11849 & 11{,}849 & $10^{9}$ & 5.4483\textsuperscript{\textdaggerdbl} & 5.7726 & 198.7 & 5.9124 & 199.4 & 5.6133 & 198.6 & \textbf{\underline{5.3876}} & 198.6 & $\mathbf{-1.11}$ & $-6.67$ & $-4.02$ \\
\addlinespace[2pt]
\midrule
\multicolumn{15}{@{}l}{\emph{Large set}}\\
usa13509 & 13{,}509 & $10^{11}$ & -- & 1.2830 & 213.9 & 1.2985 & 214.1 & 1.1939 & 213.7 & \textbf{1.1074} & 213.7 & -- & $-13.69$ & $-7.25$ \\
fma21553 & 21{,}553 & $10^{8}$ & -- & 7.9098 & 287.2 & 7.8922 & 287.0 & 7.7008 & 286.8 & \textbf{7.1350} & 286.8 & -- & $-9.80$ & $-7.35$ \\
lsb22777 & 22{,}777 & $10^{8}$ & -- & 7.7818 & 299.0 & 7.6861 & 298.2 & 7.5257 & 298.0 & \textbf{6.9678} & 298.0 & -- & $-10.46$ & $-7.41$ \\
bbz25234 & 25{,}234 & $10^{8}$ & -- & \textsc{oom} & -- & 9.5485 & 320.7 & 9.4933 & 320.3 & \textbf{8.8232} & 320.3 & -- & -- & $-7.06$ \\
boa28924 & 28{,}924 & $10^{9}$ & -- & \textsc{oom} & -- & 1.2794 & 354.6 & 1.2473 & 353.9 & \textbf{1.1620} & 353.9 & -- & -- & $-6.84$ \\
pla33810 & 33{,}810 & $10^{12}$ & -- & \textsc{oom} & -- & 1.1447 & 398.6 & 1.1194 & 398.3 & \textbf{1.0757} & 398.3 & -- & -- & $-3.90$ \\
\bottomrule
\end{tabular}
\vspace{2pt}
{\footnotesize
\textsuperscript{\textdagger}\,Proven optimal~\cite{abeledo2013} ($N\!\le\!107$). \quad
\textsuperscript{\S}\,Proven optimal~\cite{bulhoes2018} (ch150). \quad
\textsuperscript{\P}\,Best value reported by Santana et~al.~\cite{santana2022}. \quad
\textsuperscript{\textdaggerdbl}\,From~\cite{yelmewad2021}, under a different, deflated metric (cf.\
Section~\ref{sec:pdlsh-fix}); optimistic reference only.\par}
\end{threeparttable}
\end{table*}

Below $1{,}000$ nodes, GILS-RVND
holds the lowest cost on every instance, recovering the
proven optima on \texttt{berlin52}, \texttt{st70}, \texttt{kroA100}, and
\texttt{ch150}; ROAR stays within $6\%$ of reference values without attaining them. The ordering changes with scale, as the descent completes
fewer passes within $T(N)$: from $2{,}392$ nodes ROAR attains the lowest cost on
\texttt{pr2392}, \texttt{fnl4461}, \texttt{rl5934}, \texttt{pla7397}, and
\texttt{rl11849}, with a margin over GILS-RVND ranging from parity on \texttt{pr2392}
to $7.4\%$ on \texttt{pla7397}, the baseline retaining an edge only on
\texttt{pcb3038}. PDLSH holds the highest cost on every regular
instance, since distributing move evaluation across GPU threads raises the number of
moves examined per second but does not alter which moves are examined. Interestingly, although measured under the depot-inclusive metric against the optimistic references marked
$\ddagger$, ROAR falls below the reference on all eight applicable instances, LNS+2-opt
on six, and GILS-RVND on five.

\subsection{ILO versus Post-Repair Local Search}
\label{sec:ablation}
Commonly, an LNS scheme applies local search to the rebuilt tour
after each repair phase~\cite{voigt2025}, whereas ROAR applies a restricted
2-opt movement at each insertion, following the ILO paradigm. \textbf{LNS+2-opt} adopts the
common arrangement, sharing all other mechanisms while differing only in the placement of
the local search and its scope: each loop performs a standard LNS remove–reinsertion of nodes, followed by a 2-opt full neighborhood evaluation and move; this gap demonstrates the ILO impact on ROAR performance.

The ILO paradigm gives the lower mean cost on all but one regular
instance (\texttt{pr107}, below $0.01\%$), and above $1{,}000$ nodes its margin
widens with $N$, averaging $0.48\%$ from $1{,}432$ to $5{,}934$ nodes and
reaching $3.52\%$ and $4.02\%$ on \texttt{pla7397} and \texttt{rl11849}. On \texttt{rl11849}, the
$6.67\%$ gain over GILS-RVND splits into $2.76\%$ from the ruin-and-recreate
frame, while the ILO paradigm further improves ROAR performance by $4.02\%$ over its post-repair local-search counterpart.

\subsection{Large Set}
\label{sec:res-large}

On the large set (Table~\ref{tab:main}), ROAR attains the lowest cost on every
instance, by $9.8$--$13.7\%$ over GILS-RVND where it executes and by $6.0$--$14.7\%$
over PDLSH throughout. LNS+2-opt ranks second throughout, yet the margin ILO adds
over it widens with scale, from under $1\%$ below $6{,}000$ nodes to $6.8$--$7.4\%$
here (\texttt{pla33810} excepted). Notably, GILS-RVND suffers \textbf{Out-of-memory} (OOM) for instances larger than $25{,}000$ nodes, due to the auxiliary matrices for the 2-opt evaluator. PDLSH achieves inferior performance to ROAR on all instances, demonstrating that GPU acceleration cannot compensate for its search structure.

\subsection{Extreme Set}
\label{sec:extreme}
\begin{table}[t]
\centering
\begin{threeparttable}
\caption{Extreme set ($32{,}892$--$200{,}000$ nodes). Costs are divided by
the Scale factor.}
\label{tab:extreme}
\setlength{\tabcolsep}{3.5pt}
\renewcommand{\arraystretch}{1.1}
\footnotesize
\begin{tabular}{@{}l r c r r r r@{}}
\toprule
Instance & $N$ & Scale & Draft & ROAR & Red. & $t$ \\
         &     &       & (NN)  &      & (\%) & (s) \\
\midrule
xib32892      &  32{,}892 & $10^{9}$  & 1.7146 & 1.5448 & 9.91 &  389.9 \\
pba38478      &  38{,}478 & $10^{9}$  & 2.3575 & 2.1288 & 9.70 &  440.7 \\
ics39603      &  39{,}603 & $10^{9}$  & 2.3468 & 2.1601 & 7.95 &  450.9 \\
dan59296      &  59{,}296 & $10^{9}$  & 5.4106 & 4.9958 & 7.67 &  630.0 \\
pla85900      &  85{,}900 & $10^{12}$ & 6.4687 & 6.2132 & 3.95 &  871.8 \\
mona-lisa100K & 100{,}000 & $10^{11}$ & 2.9531 & 2.7087 & 8.28 & 1000.0 \\
vangogh120K   & 120{,}000 & $10^{11}$ & 4.1529 & 3.8918 & 6.29 & 1000.0 \\
venus140K     & 140{,}000 & $10^{11}$ & 4.9357 & 4.5770 & 7.27 & 1000.0 \\
pareja160K    & 160{,}000 & $10^{11}$ & 6.6229 & 6.2915 & 5.00 & 1000.0 \\
courbet180K   & 180{,}000 & $10^{11}$ & 7.6387 & 7.2287 & 5.37 & 1000.0 \\
earring200K   & 200{,}000 & $10^{11}$ & 8.7499 & 8.2167 & 6.09 & 1000.0 \\
\bottomrule
\end{tabular}
\end{threeparttable}
\end{table}

At extreme scale ROAR is assessed against its own nearest-neighbor seed, and the
reduction is stable across runs, ranging from $3.95\%$ on \texttt{pla85900} to
$9.91\%$ on \texttt{xib32892}. These latencies are upper bounds rather than optima,
since the deadline binds on every instance ($t = T(N)$) and the incumbent is still
improving at interruption. The truncation is deliberate, and reflects the operational
premise of the study: a serviceable solution within a bounded, predictable timeframe
in preference to convergence at unbounded runtime.

\subsection{Discussion}
The crossover follows from how each method spends its budget. RVND commits to a
move only after scanning an entire neighborhood, so a single pass becomes costly
as $N$ grows and few loops complete within $T(N)$; ROAR's loop, bounded by the
removal cap and the restricted local search, remains affordable and keeps
improving throughout the budget. Below $1{,}000$ nodes the descent fits the
budget and the ordering reverses. PDLSH illustrates the complementary case:
acceleration buys throughput, not a better search structure.

The memory wall is specific to GILS-RVND, not to the
literature: DBMEA~\cite{tuuszabo2020}
achieves $\mathcal{O}(1)$ time-complexity evaluation in $\mathcal{O}(N)$ space. The
contribution of Section~\ref{sec:o1} is therefore a reduction beyond that state
of the art---the proposed algorithm requires $\mathcal{O}(1)$
auxiliary space, with no preprocessing and nothing maintained as the tour changes.

\section{Conclusion}
We presented ROAR, a CPU-native ruin-and-recreate heuristic, aiming at the large and extreme-scale Traveling Repairman Problem; in addition, we propose an algorithm that evaluates each 2-opt movement in $\mathcal{O}(1)$ time in a contiguous sweep. Under a fixed operational budget, ROAR matches or surpasses GILS-RVND beyond roughly $2{,}000$ nodes, dominates both baselines on the large set---where GILS-RVND exhausts memory---and yields the first metaheuristic results for the TRP at up to $200{,}000$ nodes. Future research may focus on enhancing the performance of ROAR for small- and medium-sized instances, with the objective of achieving results that are comparable to those of the gold-standard method, GILS-RVND, on these instance classes.
\section*{Acknowledgment}
Large Language Models are used to assist in writing this paper and
generating code. All material was verified by the authors,
who take full responsibility for the paper's content.
\bibliographystyle{IEEEtran}
\bibliography{IEEEabrv, refs}

@IEEEtranBSTCTL{IEEEexample:BSTcontrol,
  CTLdash_repeated_names = "no"
}

@article{bulhoes2018,
  author  = {Bulh{\~o}es, Teobaldo and Sadykov, Ruslan and Uchoa, Eduardo},
  title   = {A branch-and-price algorithm for the minimum latency problem},
  journal = {Comput. Oper. Res.},
  volume  = {93},
  pages   = {66--78},
  month   = may,
  year    = {2018},
}

@article{yelmewad2021,
  author  = {Yelmewad, Pramod and Talawar, Basavaraj},
  title   = {Parallel deterministic local search heuristic for minimum
             latency problem},
  journal = {Cluster Comput.},
  volume  = {24},
  number  = {2},
  pages   = {969--995},
  month = jun,
  year    = {2021},
}

@article{caseau1999,
  author  = {Caseau, Yves and Laburthe, Fran{\c c}ois},
  title   = {Heuristics for large constrained vehicle routing problems},
  journal = {J. Heuristics},
  volume  = {5},
  number  = {3},
  month = oct,
  pages   = {281--303},
  year    = {1999},
}

@article{caseau2001,
  author  = {Caseau, Yves and Silverstein, Glenn and Laburthe, Fran{\c c}ois},
  title   = {Learning hybrid algorithms for vehicle routing problems},
  journal = {Theory Pract. Log. Program.},
  
  volume  = {1},
  number  = {6},
  pages   = {779--806},
  year    = {2001},

}

@article{silva2012,
  author  = {Silva, Marcos Melo and Subramanian, Anand and Vidal, Thibaut
             and Ochi, Luiz Satoru},
  title   = {A simple and effective metaheuristic for the minimum latency
             problem},
  journal = {Eur. J. Oper. Res.},
  volume  = {221},
  number  = {3},
  pages   = {513--520},
  month   = sep,
  year    = {2012},
}

@article{tuuszabo2020,
  author  = {T{\"u}{\H u}-Szab{\'o}, Boldizs{\'a}r and F{\"o}ldesi, P{\'e}ter
             and K{\'o}czy, L{\'a}szl{\'o} T.},
  title   = {An efficient evolutionary metaheuristic for the traveling
             repairman (minimum latency) problem},
  journal = {Int. J. Comput. Intell. Syst.},
  volume  = {13},
  number  = {1},
  pages   = {781--793},
  year    = {2020},
  month   = jun
}

@article{santana2022,
  author  = {Santana, {\'I}talo and Plastino, Alexandre and Rosseti, Isabel},
  title   = {Improving a state-of-the-art heuristic for the minimum latency
             problem with data mining},
  journal = {Int. Trans. Oper. Res.},
  volume  = {29},
  number  = {2},
  
  pages   = {959--986},
  year    = {2022},
  month   = mar
}

@article{voigt2025,
  author  = {Voigt, Stefan},
  title   = {A review and ranking of operators in adaptive large
             neighborhood search for vehicle routing problems},
  journal = {Eur. J. Oper. Res.},
  volume  = {322},
  number  = {2},
  pages   = {357--375},
  year    = {2025},
  month   = apr
}

@article{gendreau1992,
  author  = {Gendreau, Michel and Hertz, Alain and Laporte, Gilbert},
  title   = {New insertion and postoptimization procedures for the
             traveling salesman problem},
  journal = {Oper. Res.},
  volume  = {40},
  number  = {6},
  pages   = {1086--1094},
  year    = {1992},
  month   = dec
}

@article{vidal2014,
  author  = {Vidal, Thibaut and Crainic, Teodor Gabriel and Gendreau, Michel
             and Prins, Christian},
  title   = {A unified solution framework for multi-attribute vehicle
             routing problems},
  journal = {Eur. J. Oper. Res.},
  volume  = {234},
  number  = {3},
  pages   = {658--673},
  month   = may,
  year    = {2014}
}

@article{reinelt1991,
  author  = {Reinelt, Gerhard},
  title   = {{TSPLIB}---{A} traveling salesman problem library},
  journal = {ORSA J. Comput.},
  volume  = {3},
  number  = {4},
  pages   = {376--384},
  year    = {1991},
  month   = nov
}

@inproceedings{kaplan2005,
  author    = {Kaplan, Craig S. and Bosch, Robert},
  title     = {{TSP} art},
  booktitle = {Proc. Renaissance Banff: Math., Music, Art, Culture},
  pages     = {301--308},
  address   = {Banff, AB, Canada},
  year      = {2005},
  month     = jul
}

@article{abeledo2013,
  author  = {Abeledo, Hern{\'a}n and Fukasawa, Ricardo and Pessoa, Artur and Uchoa, Eduardo},
  title   = {The time dependent traveling salesman problem: Polyhedra and algorithm},
  journal = {Math. Program. Comput.},
  volume  = {5},
  number  = {1},
  pages   = {27--55},
  year    = {2013},
  month   = mar
}

@manual{reinelt1995,
  author       = {Reinelt, Gerhard},
  title        = {{TSPLIB} 95},
  organization = {Inst. Angew. Math., Univ. Heidelberg},
  address      = {Heidelberg, Germany},
  year         = {1995}
}

@electronic{rohe_vlsi,
  author       = {Rohe, Andr{\'e}},
  title        = {{VLSI} data sets},
  organization = {Univ. Waterloo},
  address      = {Waterloo, ON, Canada},
  url          = {https://www.math.uwaterloo.ca/tsp/vlsi/index.html},
}

@electronic{bosch_art,
  author       = {Bosch, Robert},
  title        = {{TSP} art instances},
  organization = {Univ. Waterloo},
  address      = {Waterloo, ON, Canada},
  url          = {https://www.math.uwaterloo.ca/tsp/data/art/},
}

@article{croes1958,
  author  = {Croes, Georges A.},
  title   = {A method for solving traveling-salesman problems},
  journal = {Oper. Res.},
  volume  = {6},
  number  = {6},
  pages   = {791--812},
  month   = dec,
  year    = {1958},
}

@article{rosenkrantz1977,
  author  = {Rosenkrantz, Daniel J and Stearns, Richard E and Lewis II, Philip M},
  title   = {An analysis of several heuristics for the traveling
             salesman problem},
  journal = {SIAM J. Comput.},
  volume  = {6},
  number  = {3},
  pages   = {563--581},       
  month   = sep,               
  year    = {1977},
}

@inproceedings{shaw1998,
  author    = {Shaw, Paul},
  title     = {Using constraint programming and local search methods to
               solve vehicle routing problems},
  booktitle = {Proc. 4th Int. Conf. Princ. Pract. Constraint Program.},
  pages     = {417--431},
  month     = oct,
  year      = {1998},
  address = {Pisa, Italy}
}

@article{ropke2006,
  author  = {Ropke, Stefan and Pisinger, David},
  title   = {An adaptive large neighborhood search heuristic for the
             pickup and delivery problem with time windows},
  journal = {Transp. Sci.},
  volume  = {40},
  number  = {4},
  pages   = {455--472},
  month   = nov,               
  year    = {2006},
}

@article{schrimpf2000,
  author  = {Schrimpf, Gerhard and Schneider, Johannes and
             Stamm-Wilbrandt, Hermann and Dueck, Gunter},
  title   = {Record Breaking Optimization Results Using the Ruin and
             Recreate Principle},
  journal = {J. Comput. Phys.},
  volume  = {159},
  number  = {2},
  pages   = {139--171},
  month   = apr,
  year    = {2000}
}

@inproceedings{tran2026beyond,
  title={Beyond Static Priors: Dynamic Neural Guidance for Large-Scale Ant Colony Optimization},
  author={Tran, Dat Thanh and Vu, Van Khu and Ma, Yining},
  booktitle={Proceedings of the 32nd ACM SIGKDD Conference on Knowledge Discovery and Data Mining V. 2},
  pages={4671--4682},
  year={2026}
}

@INPROCEEDINGS{11365075,
  author={Dat, Tran Thanh and Khoi, Pham Anh and Khai, Tran Quang and Dong, Do Duc and Van Khu, Vu},
  booktitle={2025 RIVF International Conference on Computing and Communication Technologies (RIVF)}, 
  title={NeuFACO: Neural Focused Ant Colony Optimization for Traveling Salesman Problem}, 
  year={2025},
  volume={},
  number={},
  pages={439-444},
  doi={10.1109/RIVF68649.2025.11365075}}
\end{document}